\documentstyle[12pt]{article}
\title{A critical study on the concept of identity in Zermelo-Fraenkel-like axioms}
\author{Aur\'elio Sartorelli$^*$ \and D\'ecio Krause$^{**}$ \and Adonai S. Sant'Anna$^*$}
\date{$^*${\it Departamento de Matem\'atica, UFPR, C.P. 019081, Curitiba, PR, 81.531-990, Brazil.}\\$^{**}${\it Departamento de Filosofia, UFSC, C.P. 476, Florian\'opolis, SC, 88040-900, Brazil.}}
\begin{document}
\maketitle

\newtheorem{definicao}{Definition}
\newtheorem{teorema}{Theorem}
\newtheorem{lema}{Lemma}
\newtheorem{corolario}{Corolary}
\newtheorem{proposicao}{Proposition}
\newtheorem{axioma}{Axiom}
\newtheorem{observacao}{Observation}

\begin{abstract}

\noindent
According to Cantor, a set is a collection into a whole of defined and separate (we shall say  {\it distinct\/}) objects. So, a natural question is ``How to treat as  `sets' collections of {\it indistinguishable\/} objects?". This is the aim of quasi-set theory, and this problem was posed as the first of present day mathematics, in the list resulting from the Congress on the Hilbert Problems in 1974. Despite this pure mathematical motivation, quasi-sets have also a strong commitment to the way quantum physics copes with elementary particles. In this paper, we discuss the axiomatics of quasi-set theory and sketch some of its applications in physics. We also show that quasi-set theory allows us a better and deeper understanding of the role of the concept of equality in mathematics.
\end{abstract}


\section{Introduction}

	It is well known that quantum theories treat elementary particles as entities which may be absolutely indistinguishable, having all their properties in common. Physicists say that they are `identical', while philosophers prefer to call them {\it indistinguishable\/}, since the standard philosophical jargon regards identical things to be the very same thing. Standard set theories like Zermelo-Fraenkel are so that the elements of a set (so as the sets themselves) obey a well defined theory of identity, according to which, roughly speaking, two objects $a$ and $b$ are always either equal or distinct; furthermore, if they are equal, then they are the very same object. If they are distinct, there exists at least one set to which one of them belongs while the another one does not (in extensional contexts, we may say that there exists at least one property which distinguishes them, a result which is usually termed Leibniz Law). In particular, the elements of a set can always be considered as {\it individuals\/} of a sort, being capable (at least ideally) of being counted, ordered or named. 

	In this sense, standard set theories cannot deal with ``genuine'' collections of indistinguishable objects. Sets are, according to the well known ``definition'' given by Cantor, ``collections of definite  and separate objects of our intuition or of our thought'' (see \cite{Cantor-55}, p. 85). The axioms of the standard set theories preserve this intuitive idea. The usual way mathematicians may consider indistinguishable things vary. The most common devices are either the consideration of permutational symmetries or group invariance \cite{Weingartner-96}. The general idea of these approaches can be summed up by the technique presented by H. Weyl \cite{Weyl-49}; in short our starting point is a set of, say, $n$ elements endowed with an equivalence relation $R$. The equivalence classes $C_{1}, \ldots, C_{k}$ were taken to play the role of collections of indistinguishable objects. But this device makes sense only if we hide the very nature of the original set as a collection of {\it distinguishable\/} objects, as given by the axioms of set theory. So, this device cannot be considered to be a right answer to the problem of finding axioms (in the sense of `set' theories) to deal with collections of {\it indistinguishable\/} objects. Similar restrictions can be made to the other techniques usually considered by mathematicians (see, for example, \cite{Sant'Anna-00a}).

	But why this kind of problem should be important? Let us recall that this is the first problem in the list proposed during the Congress on the Hilbert Problems, organized by the American Mathematical Society in 1974 (see \cite{Browder-76} p. 36). The motivation for stating this problem is of course quantum physics, which deals (whithin the scope of classical set theories!) with indistinguishable objects; as put by Yuri Manin, when presenting the just mentioned problem, 

\begin{quote}
We should consider possibilities of developing a totally new language to speak about infinity (...) I would like to point out that (...) [the usual language of set theory] is (...) an extrapolation of common-place physics, where we can distinguish things, count them, put them in some order, etc. New quantum physics has shown us models of entities with quite different behaviour. Even `sets' of photons in a looking-glass box, or of electrons in a nickel piece are much less cantorian than the `set' of grains of sand. In general, a highly probabilistic `physical infinity' looks considerably more complicated and interesting than a plain infinity of `things' (...) The twentieth century return to Middle Age scholastics taught us a lot about formalisms. Probably it is time to look outside again. Meaning is what really matters.
\end{quote}

	Of course that the problem is not only to find a way of expressing indistinguishability. Physicists do this by a trick of considering that only symmetric and antisymmetric vectors on an appropriate Hilbert space have counterpart in reality, but the interesting foundational problem is to consider indistinguishability {\it right at the start\/}, as something which is very peculiar of the objects being supposed, as in the case of quantum objects \cite{Post-63}. In this sense, we should not consider the objects first as individuals, as when we take them as elements of a set, and then to find a way of imagining them as indistinguishable entities. Indistinguishability should be a primitive concept; things should be postulated to be so that they could be indistinguishable without turnig to be identical. But, how to build a theory in which identity has no meaning (at least within some scope) but indistinguishability is meaningful if classical logic and set theories are strongly commited with Leibniz Law? This is what we present in this paper, which should be regarded as a continuation of a serie of works on the foundations of {\it quasi set theory\/} and its possible consequences in mathematics, physics, and, consequently, philosophy.

	One natural question that could be raised by the reader is: why should we use quasi-set theory in physics? According to any textbook about statistical mechanics, we know that Maxwell-Boltzmann (MB) statistics gives us the most probable distribution of $N$ {\em distinguishable\/} objects into, say, boxes with a specified number of objects in each box. We show, e.g., that the hypothesis concerning distinguishable objects is unnecessary. Usually, classical and quantum distribution functions are mathematically derived in a na\"{\i}ve fashion; but in our case an axiomatic framework is necessary if we want to show that individuality is not a necessary assumption in classical statistical mechanics. In a very interesting paper N. Huggett \cite{Huggett-99} demonstrates that the occurence of Maxwell-Boltzmann statistics in classical mechanics does not allow us to decide the metaphysical issue of molecules in a gas. In this paper we show that Maxwell-Boltzmann statistics is not committed to a metaphysical hypothesis concerning individuals.

\section{Quasi-sets}

	In trying to build a theory of collections of indistinguishable objects in the sense mentioned in the previous section, we strongly consider the motivation provided by quantum physics, with particular emphasis to Erwin Schr\"odinger's idea that the concept of identity does not make sense for elementary particles \cite{Schrodinger-52,Schrodinger-98}. In brief, this suggests that if $x$ and $y$ denote, say, electrons, it is simply meaningless to say that $x$ is identical (or different) from $y$. The lack  of meaning of such a way of speaking goes in the same direction than that one posed by W. Heinsenberg when he explains why from the point of view of modern physics the problem posed by the ancient atomists of looking for the ultimate parts of matter has no meaning. As he said (\cite{Heisenberg-89}, p. 82), 

\begin{quote}
We ask, `What does the proton consist of?' `Is the light-quantum simple,
or is it composite?' But these questions are wrongly put, since the words
{\it divide\/} or {\it consist of\/} have largely lost their meaning. It
would thus be our task to adapt our language and thought, and hence
also our scientific philosophy, to this new situation engendered by the
experiments.
\end{quote}

	Quasi-set theory was developed with the aim of providing axioms for dealing with collections of indistinguishable objects. In a certain sense, it provides an answer to the `Manin Problem' mentioned above. 

	Quasi-set theory ${\cal Q}$ is based on {\em ZFU\/}-like axioms (Zermelo-Fraenkel with {\it Urelemente\/}), but allows the existence of {\it two\/} sorts of atoms, termed respectively $m$-atoms and $M$-atoms. Two primitive unary predicates help in expressing that: $m(x)$ says that $x$ is an $m$-atom and $M(x)$ says that $x$ is an $M$-atom, where $x$ is an individual variable. The language still encompasses the binary primitive predicates $\equiv$ (indistinguishability) and $\in$ (membership), one unary functional symbol $qc$ (quasi-cardinal) and a unary predicate letter $Z$ (where $Z(x)$ says that $x$ is a {\it set\/}; these 
quasi-sets will correspond precisely to the sets of {\em ZFU\/}).  The basic idea is that the $M$-atoms  have the properties of standard {\it Urelemente\/} of {\em ZFU\/}, while the $m$-atoms may be thought of as representing elementary particles of quantum physics. Following Erwin Schrödinger, to  this last kind of entities, the concept of identity cannot be applied (\cite{Schrodinger-52}, pp. 17-8).\footnote{When we talk on `the traditional concept of identity' we mean the theory of identity as presented in standard mathematics, either in first order theories or in higher order ones (and set theory).} In quasi set theory,  this restriction is achieved by a restriction on the concept of formula: expressions like $x = y$ are not well formed if $x$ and $y$ denote $m$-atoms, despite the expression $x \equiv y$, which is read `$x$ is indistinguishable from $y$' makes sense for all the objects of the domain. The equality symbol is not primitive in our theory, but a concept of {\it extensional identity\/} is defined (see below) so that it has all the properties of standard identity of {\it ZFU}.  Then, the axiomatics permits us to distinguish between the concepts of (extensional) {\it identity\/} (being the very same object) and {\it indistinguishability\/} (agreement with respect to all the attributes), which cannot be done in classical logic and set theory.\footnote{In standard mathematics, two entities $x$ and $y$ which have the same properties  are the very same entity; this is Leibniz's Law, as already remarked.}

	A quasi-set (qset for short) $x$ is defined as something which is not a {\em Urelemente\/}. A qset $x$  may have a cardinal (termed its {\it quasi cardinal\/}, and denoted by $qc(x)$) but, in general, it has not an ordinal, since there are quasi sets which cannot be ordered 
(since their elements are indistinguishable $m$-atoms, expressed by the relation $\equiv$).  The concept of quasi cardinal is taken as primitive, since it cannot be defined by usual means.  This fits the idea that quantum particles cannot be ordered or counted, but only aggregated in certain amounts. Notwithstanding, due to the concept of quasi cardinal, there is a sense (as in ortodox quantum physics) in saying that there may exist a certain quantity of $m$-atoms obeying certain conditions, despite they cannot be named or labeled.  

	The primitive relation of indistinguishability ($\equiv$) is postulated to be reflexive, symmetric and transitive, but in order to differenciate it from identity as ascribed by the traditional (first-order) theory of identity, the substitutivity axiom does not hold. Even so, it should be interesting that such a relation, which holds in the whole domain, turns to be the standard identity (here represented by the extensional identity defined below) when the objects under consideration are not $m$-atoms. Then a concept of {\it extensional identity\/} fits the idea of classical identity. The first definitions and axioms are the following:

\begin{definicao}\label{defqset}
\hfill{ }
\begin{enumerate}

\item $Q(x) := \neg (m(x) \vee M(x))$. We read $Q(x)$ as ``$x$ is a quasi-set'' or ``$x$ is a qset''.

\item $P(x) := Q(x) \wedge \forall y (y \in x \Rightarrow m(y)) \wedge \forall y \forall z (y \in x \wedge z \in x \Rightarrow y \equiv z)$. In this case we say that $x$ is a {\em pure qset\/}.

\item $D(x) := M(x) \vee Z(x)$. These are the `(classical) things', to use Zermelo's original terminology. We read $D(x)$ as ``$x$ is a {\em Dinge/}''.

\item $E(x)  := Q(x) \wedge \forall y (y \in x \Rightarrow Q(y))$.

\item $x =_{E} y := (Q(x) \wedge Q(y) \wedge \forall z ( z \in x 
\Leftrightarrow z \in y )) \vee (M(x) \wedge M(y) \wedge \forall _Q z (x\in z \Leftrightarrow y\in z))$. In this case we say that $x$ and $y$ are extensionaly identical. The symbol ``$=_E$'' is called extensional identity.

\item $x \subseteq y := \forall z (z \in x \Rightarrow z \in y)$.

\end{enumerate}

\end{definicao}

	The first axioms of ${\cal Q}$ are: 

\vspace{2mm}
\noindent 
{\bf (Q1)} $\forall x (x \equiv x)$

\vspace{3mm}
\noindent 
{\bf (Q2)} $\forall x \forall y (x \equiv y \Rightarrow y \equiv x)$

\vspace{3mm}
\noindent 
{\bf (Q3)} $\forall x \forall y \forall z (x \equiv y \wedge y
\equiv z \Rightarrow x \equiv z)$

\vspace{3mm}
\noindent 
{\bf (Q4)} $\forall x \forall y (x =_{E} y \Rightarrow (A(x,x) \Rightarrow
A(x,y)))$, with the usual syntactic restrictions.

\begin{teorema}
Whether $Q(x)$ or $M(x)$, then $x =_{E} x$.
\end{teorema}

\begin{description}
\item [Proof:] If $Q(x)$, since $\forall z (z \in x \Leftrightarrow z \in x)$, then $x =_{E} x$ by the definition of extensional identity. If $M(x)$,
then since $x \equiv x$ by {\bf Q1}, it follows that $x =_{E} x$. $\Box$
\end{description}

\begin{corolario}
The relation of extensional equality has all the properties 
of classical equality.
\end{corolario}

\begin{description}
\item[Proof:] Straightforward, if we take into account the above theorem and {\bf Q4}. $\Box$
\end{description}

\vspace{3mm}
\noindent
{\bf (Q5)} Nothing is at the same time an $m$-atom and an $M$-atom:
$$\forall x (\neg (m(x) \wedge M(x)))$$

\begin{teorema}\label{notm}
Wheter $Q(x)$ or $M(x)$, then $\neg m(x)$.
\end{teorema}

\begin{description}
\item [Proof:] If $Q(x)$, then $\neg m(x)$ by the definition of qset.
If $M(x)$, then $\neg m(x)$ by {\bf Q5}. $\Box$
\end{description}

\vspace{3mm}
\noindent
{\bf (Q6)} The atoms are empty: 
$$\forall x \forall y (x \in y \Rightarrow Q(y))$$

\vspace{3mm}
\noindent
{\bf (Q7)} Every set is a qset:
$$\forall x (Z(x) \Rightarrow Q(x))$$

\vspace{3mm}
\noindent
{\bf (Q8)} Qsets whose elements are `classical things' are
sets and conversely:
$$\forall_{Q} x (\forall y (y \in x \Rightarrow D(y)) 
\Leftrightarrow Z(x))$$

\noindent
What is the meaning of {\bf Q8}? Our intention is to characterize the {\it sets\/} in ${\cal Q}$ so that they can be identifyed with the sets of {\em ZFU\/}. This is supposed to be achieved if they were taken to be those qsets whose transitive closure (this concept can be defined in the usual sense) does not contain $m$-atoms.  The `$\Rightarrow$-part' of {\bf Q8} gives half of the answer: if all the elements of $x$ are {\it Dinge\/} (either sets of $M$-atoms), then $x$ is a set. Concerning the converse, it is not enough to postulate that  no element of a set is an $m$-atom, since  it may be that the elements of its elements have $m$-atoms as elements and so on. The problem can be satisfactorily solved  if we have $Z(x) \Rightarrow \forall y (y \in x \Rightarrow D(y))$, which is precisely the `$\Leftarrow$-part' of {\bf Q8}. 

\vspace{0.3cm}
\noindent
{\bf (Q9)}  
$$\forall x (m(x) \wedge x \equiv  y \Rightarrow m(y)) \wedge
\forall x \forall y (x =_{E} y \wedge M(x) \Rightarrow M(y))$$
$$\wedge \forall x \forall y (x =_{E} y \wedge Z(x) \Rightarrow Z(y))$$

\noindent
{\bf (Q10)} The empty qset: there exists a qset denoted by `$\emptyset$', which does not have elements: 
$$\exists_{Q} x \forall y (\neg (y \in x))$$ 

\vspace{0.3cm}
\begin{teorema}
The empty qset is a set.
\end{teorema}

\begin{description}
\item [Proof:] Take $x =_{E} \emptyset$. Since $y \in x$ is false by {\bf Q10}, then the antecedent of $\forall y (y \in x \Rightarrow D(x))$ is true, hence $Z(x)$ by {\bf Q8}. $\Box$
\end{description}

\vspace{0.3cm}
\noindent
{\bf (Q11)} Indistinguishable {\it Dinge\/} are extensionally identicals:
$$\forall_{D} x \forall_{D} y (x \equiv y \Rightarrow x =_{E} y)$$

\vspace{0.3cm}
\noindent
{\bf (Q12)} This is the qset-theoretical version of the weak-pair axiom. For all  $x$ and $y$, there exists a qset whose elements are the indistinguishable from either $x$ or $y$: $$\forall x \forall y \exists_{Q} z \forall t (t \in z \Leftrightarrow t \equiv x \vee t \equiv y)$$

\noindent
We denote this qset by $[x, y]$ and $\{x, y \}$ when $x$ and $y$ are 
{\it Dinge\/}, as usual.

	As we will see below after giving the idea of the quantity of elements of a qset (by means of the primitive concept of quasi-cardinal), the quasi-cardinal of $[x]$ may be different from 1, where $[x]:= [x,x]$. 

\vspace{0.3cm}
\noindent
{\bf (Q13)} The Separation Schema: by considering the usual syntactical restrictions on the formula $A(t)$, the following is an axiom: 
$$\forall_{Q} x \exists_{Q} y \forall t (t \in y \Leftrightarrow t \in x \wedge A(t))$$

\vspace{3mm}
This qset is written
$[t \in x : A(t)]$ (we may use \{ and \} when such a qset is a set).

\vspace{0.3cm}
\noindent
{\bf (Q14)} Union 
$$\forall_{Q} x (E(x) \Rightarrow \exists_{Q} y
(\forall z (z \in y) \Leftrightarrow \exists t (z \in t \wedge t \in x)))$$

\vspace{3mm}
This qset is denoted by  $\bigcup_{t \in x} t$ \, (we also
use $x \cup y$ as usual).

\vspace{0.3cm}
\noindent
{\bf (Q15)} Power-qset 
$$\forall_{Q} x \exists_{Q} y \forall t
(t\in y \Leftrightarrow t \subseteq x)$$ 

\vspace{3mm}
According to the standard notation, we write ${\cal P}(x)$ for this qset.

\begin{definicao}
\hfill{ }
\begin{enumerate}

\item $\langle x, y \rangle := [[x],[x,y]]$ 

\item  $x \times y := [\langle
z, u \rangle \in {\cal P}{\cal P}(x \cup y) : z \in x \wedge u
\in y]$

\item The concepts of intersection and difference of qsets are
defined in the usual way so that $t \in x \cap y$ iff $ t \in x \wedge t \in y$ and $t \in x - y$ iff $t \in x \wedge t \notin y$. This last concept will be mentioned again later. It is worth to note that the symbol `$\notin$' has its usual meaning as in set theory.

\end{enumerate}

\end{definicao}

	We remark that $\langle x, y \rangle$ is a kind of `generalized ordered pair', since the first element is the qset of all indistinguishable from $x$, while the second is the qset of all indistinguishable from $y$. We call it the `weak pair'. The collection $[x]$ is termed a `weak singleton' of $x$.

\vspace{0.3cm}
\noindent
{\bf (Q16)} Infinity: $$\exists_{Q} x (\emptyset \in x \wedge \forall y (y \in x \wedge Q(y) \Rightarrow y \cup [y] \in x))$$.

\vspace{0.3cm}
\noindent
{\bf (Q17)} Regularity: (Qsets are well-founded):
$$\forall_{Q} x (E(x) \wedge x \neq_{E} \emptyset \Rightarrow \exists_{Q} y (y \in x \wedge y \cap x =_{E} \emptyset))$$

\subsection{Relations} 

\begin{definicao}
A qset $w$ is a {\it relation} if it satisfies the following predicate $R$: 
$$R(w) := Q(w) \wedge \forall z (z \in w \Rightarrow \exists u
\exists v (u \in x \wedge v \in y \wedge z =_{E} \langle u, v
\rangle))$$
\end{definicao}

\begin{teorema}
No partial, total or strict order relation can be defined on a pure qset whose elements are indistinguishable from one another.
\end{teorema}

\begin{description}
\item [Proof:] (Sketch) Partial and total orders require antisymmetry, and this property cannot be stated without identity. Asymmetry also cannot be supposed, for, if $x \equiv y$, then for every $R$ such that $\langle x, y \rangle \in R$, and there it follows that $\langle x, y \rangle =_{E} [[x]] =_{E} \langle y, x \rangle \in R$; so, $xRy$ entails $yRx$. $\Box$
\end{description}

\begin{teorema}
For every formula $A$ of {\em ZFU\/} let $A^{q}$ be its translation to the language of ${\cal Q}$, then  $\vdash_{ZFU} A$ iff $\vdash_{\cal Q} A^{q}$.
\end{teorema}

\begin{description}
\item [Proof:] The theory ${\cal Q}$ encompasses a `classical' counterpart which can be defined as follows: let $A$ be a formula of the language of {\em ZFU\/} (which we may admit has an unary predicate $S$ which stands for `sets'. Then, call $A^{q}$ its translation to ${\cal Q}$, defined as follows, where $S(x)$ means that $x$ is a set (in ZFU):

\begin{enumerate}
\item If $A$ is $S(x)$, then $A^{q}$ is $Z(x)$
\item If $A$ is $x = y$, then $A^{q}$ is $((M(x) \wedge M(y)) 
\vee (Z(y) \wedge Z(y)) \wedge x =_{E} y)$
\item If $A$ is $x \in y$, then  $A^{q}$ is $((M(x) \vee
Z(x)) \wedge Z(y)) \wedge x \in y$
\item If $A$ is $\neg B$, then $A^{q}$ is $\neg B^{q}$
\item If $A$ is $B \vee C$, then $A^{q}$ is $B^{q} \vee C^{q}$
\item If $A$ is $\forall x B$, then $A^{q}$ is $\forall x
(M(x) \vee Z(x) \Rightarrow B)$
\end{enumerate}

	Then it is easy to see that the translations of the axioms of
{\em ZFU\/} are theorems of  ${\cal Q}$. So, if ${\cal Q}$ 
is consistent, so is {\em ZFU\/} (see \cite{daCosta-99}). $\Box$
\end{description}

	The result above shows that there is a {\it copy\/} of {\em ZFU\/} in ${\cal Q}$. In this  `copy', we may define the following concepts: $Cd(x)$ for `$x$ is a cardinal';  $card(x)$ denotes `the cardinal of $x$, and $Fin(x)$ says that `$x$ is a finite quasi-set' (that is, $qc(x)$ is a natural number). 

	By considering these concepts, we may present the axioms for quasi-cardinals:\\

\noindent
{\bf (Q18)} Every object which is not a qset (that is, every 
{\it Urelement\/}) has quasi-cardinal zero:
$$\forall x (\neg Q(x) \Rightarrow qc(x) =_{E} 0)$$

\noindent
{\bf (Q19)} The quasi-cardinal of a qset is a cardinal (defined in the `classical part' of the theory and coincides with the cardinal itself when this qset is a set:
$$\forall_{Q} x \exists ! y (Cd(y) \wedge y =_{E} qc(x) \wedge (Z(x) \Rightarrow y =_{E} card(x)))$$

\noindent
{\bf (Q20)} Every non-empty qset has a non null quasi-cardinal:
$$\forall_{Q} x (x \neq_{E} \emptyset \Rightarrow qc(x) \neq_{E} 0)$$

\noindent
{\bf (Q21)} $\forall_{Q} x (qc(x) =_{E} \alpha \Rightarrow \forall \beta (\beta \leq_{E} \alpha \Rightarrow \exists_{Q} y (y \subseteq x \wedge qc(y) =_{E} \beta))$

\noindent
{\bf (Q22)} $\forall_{Q} x \forall_{Q} y \forall t (y \subseteq x \rightarrow qc(y) \leq_{E} qc(x))$

\noindent
{\bf (Q23)} $\forall_{Q} x \forall_{Q} y (Fin(x) \wedge x \subset y \Rightarrow qc(x) < qc(y))$ 

\noindent
{\bf (Q24)} $\forall_{Q} x \forall_{Q} y (\forall w (w \notin x \vee w \notin y) \Rightarrow qc(x \cup y) =_{E} qc(x) + qc(y))$\\

	In the next axiom, $2^{qc(x)}$ denotes (intuitively) the quantity of
subquasi-sets of $x$. Then,\\

\noindent
{\bf (Q25)} $\forall_{Q} x (qc({\cal P}(x)) =_{E} 2^{qc(x)})$\\

	Axiom {\bf Q25} raises an interesting discussion we mention below. But first we need the concept of {\it quasi-function\/}.

\subsection{Quasi-functions}

	Standard functions could not distinguish between arguments and values. So, we pose:

\begin{definicao}
If $x$ and $y$ are qsets and $R$ is the predicate for `relation' defined above, we say that $f$ is a {\em quasi-function\/} (qfunction) if it satisfies the following predicate: 

$$QF(f) := R(f) \wedge \forall u (u \in x \Rightarrow \exists v (v \in y \wedge
\langle u, v \rangle \in f)) \wedge $$
$$ \forall u \forall u' \forall v \forall v' (\langle u, v \rangle \in f \wedge \langle u', v' \rangle \in f \wedge u \equiv u' \Rightarrow v \equiv v')$$
\end{definicao}

\noindent
$f$ is a $q$-injection if $f$ is a $q$-function from $x$ to $y$ and satisfies the additional condition: 

$\forall u \forall u' \forall v \forall
v' (\langle u, v \rangle \in f \wedge \langle u', v' \rangle \in f \wedge v \equiv v' \Rightarrow u \equiv u') \wedge  qc(Dom(f)) \leq_{E} qc(Rang(f))$

\noindent
$f$ is a $q$-surjection if it is a function from $x$ to $y$ such
that 

$\forall v (v \in y \Rightarrow
\exists u (u \in x \wedge \langle u, v \rangle \in f)) \wedge
qc(Dem(f)) \geq_{E} qc(Rang(f)).$ 

\noindent
A function $f$ which is both a $q$-injection and a $q$-surjection is said to be a $q$-bijection.

	In this case, $qc(Dom(f)) =_{E} qc(Rang(f))$.

\subsection{How many subquasi-sets?} 

	Now we can turn to the discussion involving the axiom {\bf Q25}. Since the concept of identity has no meaning for $m$-atoms, how can we ensure that a qset $x$ such that $qc(x) =_{E} \alpha$ has precisely $2^{\alpha}$ subqsets? In standard set theories (so as in the `classical part' of ${\cal Q}$, that is, in considering those qsets which fit the sets of {\em ZFU\/}), as it is well known, if $card(x)$ denotes the cardinal of $x$, then by the definition of exponentiation of cardinals, $2^{card(x)}$ is defined to be the cardinal of the set $^{x}2$, which is the set of all functions from $x$ to the Boolean algebra $2 = \{0, 1\}$ (see \cite{Enderton-77}). In  ${\cal Q}$ this definition doesn't work. Let us explain why.

	Suppose that $\alpha$ is the quasi-cardinal of $x$, which is a cardinal, by the axiom {\bf Q19}.  This axiom says that every qset has a unique quasi-cardinal which is a cardinal (defined in the `classical part' of the theory), and if the qset is in particular a set (in ${\cal Q}$), then this quasi-cardinal is its cardinal {\em stricto sensu\/}. So, every quasi-cardinal is a cardinal and the above expression `there is a unique' makes sense. Furthermore, from the fact that $\emptyset$ is a set, it follows that its quasi-cardinal is 0. Then we may write 

\begin{equation}\label{2qcx}
2^{qc(x)} := qc(^{\alpha}2)
\end{equation}

\noindent
and then, since $\alpha$ is a cardinal and both $\alpha$ and $2$ are {\it Q-sets\/}, we have

\begin{equation}\label{eq3}
2^{qc(x)} := card(^{\alpha}2)
\end{equation}

	So, we may take the cardinal of the qset $^{\alpha}2$ in its usual sense to mean $2^{qc(x)}$. Then, (\ref{eq3}) gives meaning to the axiom {\bf Q25}, since it explains what does $2^{qc(x)}$ mean: it is the cardinal of the {\it set\/} of all the applications from $\alpha$ (the quasi-cardinal of $x$) in $2$. By considering this, the axiom may be written as follows, where $x$ is a qset and $\alpha$ is its quasi-cardinal:\\

\noindent
{\bf Axiom Q25} (Alternative Form) \,  
$$\forall_{Q} x (qc({\cal P}(x))) =_{E} card(^{\alpha}2).$$

	We remark that the second member of the equality has a precise meaning in ${\cal Q}$ , since both $\alpha$ and $2$ act as in classical set theories, as remarked above. This characterization  allows us to avoid another problem, which could be thought to be derived in the quasi-set theory. To explain it, we recall that in standard set theories we can prove that ${\cal P}(x)$ is equinumerous with $^{x}2$ by defining a one-one function $f : {\cal P}(x) \to ^x 2$  as follows: for every $y \subseteq x$, let $f(y)$ be the characteristic function of $y$, namely, the function $\chi_{y} : x \to 2$ defined by

\begin{equation}\label{chacf}
\chi_{y}(t) :=  \left\{
\begin{array}{ll}
1 & \mbox{if $t \in y$}\\
0 & \mbox{if $t \in x - y$}
\end{array}
\right.
\end{equation}

	Then any function $h \in \, ^{x}2$ belongs to the range of $f$ since
$$h = f(\{t \in x : h(t) = 1\}).$$
 
	Suppose now that $x$ is a qset such that $qc(x)$ is the natural number $n$ and that all elements of $x$ are indistinguishable one each other (the natural numbers are defined in ${\cal Q}$ in the usual way, just in the model of {\em ZFU} we have defined in ${\cal Q}$).\footnote{For all we need, it is sufficient to consider finite qsets (this definition is also standard).} In this case, we cannot define the characteristic quasi-function
$\chi^{q}_{y}$ for $y \subseteq x$, since, for instance, if $\chi^{q}_{y}(t) =_{E} 1$ for $t \in y$, then $\chi^{q}_{y}(w) =_{E} 1$ as well for every $w \in x$, independently if either $w$ belongs to $y$ or not. This is due to the definition of the quasi-functions given above, since for every quasi-function $f$,

$$\langle a, b \rangle \in f \wedge \langle c, d \rangle \in f \wedge a \equiv c \Rightarrow b \equiv d.$$

	In other words, if the image of a certain $t$ by the quasi-function $f$ is $1$, then the image of every element indistinguishable from $t$ will be $1$ as well. So, ${\cal Q}$ distinguishes only between {\it two\/} quasi-functions from $x$ to $2$, namely, that one which associate $1$ to all elements of $x$ and that one which associate $0$ to all of them. This is the motive why we used $qc(^{\alpha}2)$ to mean $2^{qc(x)}$, since both $\alpha$ and $2$ may be viewed as  {\it sets\/} (in the standard sense). If we had used $^{x}2$ instead, we would be unable to distinguish among certain quasi-functions, so complicating the meaning of {\bf Q25}, since we could have no manner of counting the number of subquasi-sets of a qset. But, by using $^{\alpha}2$, since both $\alpha$ and $2$ behave `classically', we keep {\bf Q25} with its usual meaning. 

	From these considerations, we may conclude that when $x$ is a qset whose elements are indistinguishable $m$-atoms, we cannot prove {\it within\/} ${\cal Q}$ that if $qc(x) =_{E} n$, so we cannot assert that $x$ has $2^{n}$ subquasi-sets. Since this is precisely what {\bf Q25} intuitively says, we may affirm that this axiom cannot be proved from the remaining axioms of ${\cal Q}$. But, since it holds for particular qsets, namely, to those which are {\it sets\/}, it cannot be disproved as well. In order to state that {\bf Q25} cannot be disproved, consider the {\it sets\/} in ${\cal Q}$; since they behave as classical sets, we can prove that what {\bf Q25} asserts is true. Now it suffice to take a qset whose elements are indistinguishable $m$-atoms and such that $qc(x) = \alpha$.

\subsection{The `weak' extensionality}

	The absence of a theory of identity for the $m$-atoms, due to the lack of meaning of speaking about either the identity or the difference of $m$-atoms, causes the necessity of a modification in the Axiom of Extensionality, which here does not hold as in standard set theories. In order to do so, let us introduce the following definition:

\begin{definicao}\label{sim}
For all non empty quasi-sets  $x$ and $y$,  

$Sim(x,y) := \forall z \forall t (z \in x \wedge t \in y \Rightarrow z \equiv t)$
In this case we say that $x$ and $y$ are {\em similars\/}.

\vspace{2mm}
$QSim(x,y) := Sim(x,y) \wedge qcard(x) =_{E} qcard(y)$. That is, $x$ and
$y$ are $Q$-{\em similar\/} iff they are similar and have the same quasi-cardinality.
\end{definicao}

\noindent
{\bf (Q26)} Weak Extensionality: Qsets which have the same quantity of elements of the same sort are indistinguishable. In symbols,

$$\forall_{Q} x \forall_{Q} y ((\forall z (z \in x/{\equiv} \Rightarrow 
\exists t (t \in y/{\equiv} \wedge \wedge QSim(z,t))))$$

$$\wedge \forall t (t \in y/{\equiv} \Rightarrow \exists z
(z \in x/{\equiv} \wedge \wedge QSim(t,z))) \Rightarrow x \equiv y)$$

	Axiom {\bf Q26} allows us to remark another point about {\bf Q25}. As in standard set theories, if $card(x) =_E n$, then are there exactly $n$ subsets of $x$ which are singletons? If not, how can we make sense to the idea that if $qc(x) =_{E} n$, then $x$ has $n$ elements? We recall that the main motivation of ${\cal Q}$ is the way quantum mechanics deals with elementary particles. In this theory, despite there is a sense in saying that, say, there are $k$ electrons in a certain level of a certain atom, there is no way of counting them or of distinguishing among them (see \cite{vanFraassen-91}, Chap. 12).

	If $x$ is a qset whose elements are indistinguishable from one another as above (let us suppose again that $qc(x) =_{E} n$, which suffices for our purposes), then the  singletons $y \subseteq x$ are indistinguishable one each other, as results from the weak extensionality axiom {\bf Q26}. So, all the singletons (in the intuitive sense) seem to fall in just one qset. But it should be recalled that these `singletons' (subqsets whose quasi-cardinality is $1$) are not {\it identical\/} (that is, they are not {\it the same object\/}), but they are indistinguishable in a precise sense (given by {\bf Q26}). In other words, despite the theory cannot distinguish among them, we cannot state neither that they are the same qsets nor that their elements are identical. So, it is consistent with ${\cal Q}$ to suppose that if $qc(x) = \alpha$, then $x$ has precisely $\alpha$ `singletons'. So, due to {\bf Q25}, the theory does not forbid the existence of such singletons, despite in ${\cal Q}$ we cannot prove that they exist as `distinct' entities, and hence we may reason in ${\cal Q}$ as physicists do when informaly dealing with a certain number of indistinguishable elementary particles.

	By means of {\bf Q26} it is easy to prove the following theorem:   

\begin{teorema}\label{oquefaltava}
\hfill{ }
\begin{enumerate}
\item $x =_{E} \emptyset \wedge y =_{E} \emptyset \Rightarrow x \equiv y$
\item $\forall_{Q} x \forall_{Q} y (Sim(x,y) \wedge qc(x) =_{E} qc(y) \Rightarrow 
x \equiv y)$
\item $\forall_{Q} x \forall_{Q} y (\forall z (z \in x
\Leftrightarrow z \in y) \Rightarrow x \equiv y)$
\item  $x \equiv y \wedge qc([x]) =_{E} qc([y]) \Leftrightarrow
[x] \equiv [y]$
\end{enumerate}
\end{teorema}

\subsection{Replacement axioms}

	To keep ${\cal Q}$ with a structure similar to ZFU, we may state the Replacement Axioms as it follows:

	If $A(x,y)$ is a formula in which $x$ and $y$ are free variables, we say that $A(x,y)$ defines a $y-(qfunctional)$ condition on the quasi-set $t$ if $\forall w (w \in t \Rightarrow \exists s A(w,s) \wedge \forall w \forall w' (w \in t \wedge w' \in t \Rightarrow \forall s \forall s' (A(w,s) \wedge A(w',s') \wedge w \equiv w' \Rightarrow s \equiv s'))$ (this is abbreviated by $\forall x \exists ! y A(x,y)$). Then, we have:\\

\noindent
{\bf (Q27)} Replacement

$\forall x \exists ! y A(x,y) \Rightarrow \forall_{Q} u \exists_{Q} v (\forall z
(z \in v \Rightarrow \exists w (w \in u \wedge A(w,z)))$

\subsection{The concept of strong singleton}

\begin{definicao}
A {\em strong singleton\/} of $x$ is a quasi-set $x'$ which satisfies the following property: 
$$x' \subseteq [x] \wedge qc(x') =_{E} 1$$
\end{definicao}

	In words, a strong singleton of $x$ is a qset whose only element is an indistinguishable from $x$. In standard set theories, this qset is of course the singleton whose only element is $x$ itself, but here $x$ may be an $m$-atom, and in this case there is no way of speaking of something {\it being\/} $x$.  Even so, we can prove that such a qset exists:

\begin{teorema}
For all $x$, there exists a strong singleton of $x$.
\end{teorema}

\begin{description}
\item[Proof:] The qset $[x]$ exists by the weak pair axiom. Since $x \in [x]$ (recall that $\equiv$ is reflexive), we have $qc([x]) \geq_{E} 1$ by {\bf Q20}. But, by {\bf Q21}, there exists a subqset of $[x]$ which has quasi-cardinal $1$. Take this qset to be $x'$.
$\Box$
\end{description}

\begin{teorema}
All the strong singletons of $x$ are indistinguishable.
\end{teorema}

\begin{description}
\item [Proof:] Immediate consequence of {\bf Q26}, since all of them have the same quasi-cardinality $1$ and their elements are indistinguishable by definition. $\Box$
\end{description}

	The important remark is that, as we shall see, we cannot prove that the strong singletons of $x$ are extensionaly identicals. In ${\cal Q}$, the concept of difference of qsets is introduced in the usual way: $x - y$ is the qset whose elements are the elements of $x$ which do not belong to $y$. But, in what respects indistinguishable $m$-atoms, we cannot give ostensive definitions, say by puting the finger over an $m$-atom and saying `That is Peter'. Even so, as in quantum physics, we may reason as if a certain element does belong to the qset or not; the excluded middle law remains valid, even if we cannot verify what case holds.\footnote{As in standard mathematics, let us remark, there is no effective procedure  for proving the valitity of the Law of the Excluded Middle, given inputs $x$ and $y$, to decide either $x \in y$ or $x \notin y$.} This idea fits what happens with the electrons in an atom; in general we know how many electrons there are, and we can say that some of them {\it are\/} in that atom, but we cannot tell what are the particular electrons which are in: this question simply loose its usual meaning.

	Coming back to the meaning of $x-y$, we will show that the quasi-cardinal of $x - y$ is, as expected, $qc(x) - qc(y)$. 

\begin{teorema}
For all qsets $x$ and $y$, if $y \subseteq x$, then $qc(x - y) =_{E} qc(x) - qc(y)$.
\end{teorema}

\begin{description}
\item [Proof:] By definition, $t \in x - y$ iff $t \in x \wedge t \notin y$. Then $(x - y) \cap y =_{E} \emptyset$. Hence, by {\bf Q24}, $qc((x - y) \cup y) =_{E} qc(x - y) + qc(y)$ (let us call this expression (i)). But, since $y \subseteq x$,  $(x - y) \cup y =_{E} x$ and so, in order (i) to be true,  $qc(x - y) =_{E} qc(x) - qc(y)$.$\Box$
\end{description}

	The next result may be viewed as a quasi-set-theoretical version of the indistinguishability Postulate used in quantum physics. Roughly speaking, it says that permutations of indistinguishable quanta are not observable, and constitute one of the most basic metaphysical assumptions which underly quantum mechanics \cite{Krause-95}. In order to state and prove this result, we introduce a definiton. 

\begin{definicao}\label{star}
\begin{enumerate}
\item Let $x$ be a qset such that $E(x)$, that is (according to Definition 1), its elements are also qsets. Then, 
$$\bigcap_{t \in x} t := [z \in \bigcup_{t \in x} t : \forall s (s \in x \Rightarrow z \in s)]$$
\item If $m(u)$,\footnote{This is of course the most interesting case. The generalization of this definition to $M$-atoms and sets, so as the results that follow, is immediate, but in this case the results coincide with the analogous (and sometimes trivial) situations in standard set theories. The case of $m$-atoms is really that one which makes the difference.} then $S_{u} := [s \in {\cal P}([u]) : u \in s]$
\item $u^{*} := \bigcap_{t \in S_{u}} t$
\end{enumerate}
\end{definicao}

\begin{lema}\label{Lemma} 
If $m(u)$, then:
\begin{enumerate}
\item $u \in \bigcup_{t \in S_{u}} t$
\item $\forall s (s \in S_{u} \Rightarrow u \in s)$
\item $z \in u^{*}$ iff $z \in \bigcup_{t \in S_{u}} t \wedge \forall s (s \in S_{u} \Rightarrow z \in s)$
\item $u \in u^{*}$
\item $u^{*} \subseteq [u]$
\item If $s \in S_{u}$, then $u^{*} \subseteq s$
\end{enumerate}
\end{lema}

\begin{description}
\item [Proof:] (1) $z \in \bigcup_{t \in S_{u}} t$ iff $\exists t (t \in S_{u} \wedge z \in t)$.  Therefore, by the above definition, $z \in \bigcup_{t \in S_{u}} t$ iff $\exists t (t \in {\cal P}([u]) \wedge u \in t \wedge z \in t)$. But since $[u] \in {\cal P}([u])$ and $u \in [u]$, it follows that $u \in \bigcup_{t \in S_{u}} t$. (2) $\forall s (s \in S_{u} \Leftrightarrow s \in {\cal P}([u]) \wedge u \in s)$. Therefore, $\forall s (s \in S_{u} \Rightarrow u \in s)$. (c) Immediate consequence of the above definition. (4) Immediate consequence of (1)-(3) above. (5) Suppose that $z \in u^{*}$. By (3), we have $\forall s (s \in S_{u} \Rightarrow z \in s)$. But since $[u] \in S_{u}$, it results that $z \in [u]$. (6) If $z \in u^{*}$, then, as before, $\forall s (s \in S_{u} \Rightarrow z \in s)$. But, by hyphotesis, $s \in S_{u}$; so, $z \in s$. $\Box$
\end{description}

\begin{lema}\label{lemma}
If $u$ is an $m$-atom and $z$ is a qset, then if $z \subseteq u^{*}$ and $qc(z) =_{E} 1$, it results that either $u \in u^{*} - z$ or $qc(u^{*}) =_{E} 1$.
\end{lema}

\begin{description}
\item [Proof:] Suppose $u \notin u^{*} - z$. Since $u \in u^{*}$, it follows that $u \in z$. But $z \subseteq u^{*} \subseteq [u]$, therefore $z \in S_{u}$. But, by item (6) of the above Lemma, $u^{*} \subseteq z$. By hypothesis, $z \subseteq u^{*}$, hence $u^{*} =_{E} z$, and so $qc(u^{*}) =_{E} qc(z) =_{E} 1$. $\Box$
\end{description}

\begin{teorema}\label{theorem}
For every $u$, $qc(u^{*}) =_{E} 1$.
\end{teorema}

\begin{description}
\item [Proof:] By item (4) of Lemma (\ref{Lemma}), $u^{*} \not=_{E} \emptyset$. So, by {\bf Q20}, $qc(u^{*}) \not=_{E} 0$, hence $qc(u^{*}) \geq_{E} 1$. We shall show that the equality holds. Suppose that $qc(u^*) >_{E} 1$. Then, by {\bf Q21}, there exists a qset $w \subseteq u^*$ such that $qc(w) =_{E} 1$. So, by Lemma (\ref{lemma}), $u \in u^{*} - w$. But $u^{*} - w \subseteq [u]$, since $u^{*} \subseteq [u]$. Therefore, $u^* - w \in S_u$. By Lemma (\ref{Lemma}), item (6), $u^{*} \subseteq u^{*} - w$. But since $u^{*} - w \subseteq u^{*}$, it follows that $u^{*} =_{E} u^{*} - w$. Again by {\bf Q20}, $w \not=_{E} \emptyset$ since $qc(w) =_{E} 1$. Then let be $t \in w$. So, $t \in u^{*}$ since $w \subseteq u^{*}$, hence $t \in u^{*} - w$ (once $u^{*} =_{E} u^{*} - w$). Then $t \notin w$, a contradiction. $\Box$
\end{description}

\begin{lema}\label{pinga}
For all $m$-atoms $u$ and $v$, if $u \equiv v$, then $u^{*} \equiv v^{*}$. Furthermore, if $u \in w$, then $u^{*} \subseteq w$ for any qset $w$.
\end{lema}

\begin{description}
\item [Proof:] By Lemma (1), item (5), $u^*\subseteq [u]$ and $v^*\subseteq [v]$; if $u\equiv v$ then $Sim(u^{*}, v^{*})$ (see Definition (\ref{sim})). But, by Theorem (\ref{theorem}), $qc(u^{*}) =_{E} 1$ and $qc(v^{*}) =_{E} 1$ and then, by theorem (\ref{oquefaltava}), item (2), $u^{*} \equiv v^{*}$. The last part can be proven by noting that if $u \in w$, then $u \in w \cap [u]$, so as $w \cap [u] \subseteq [u]$, therefore $w \cap [u] \in S_{u}$. Then, by Lemma (\ref{Lemma}), item (6), $u^* \subseteq w \cap [u]$ and so $u^{*} \subseteq w$. $\Box$
\end{description}

	These last results show that $u^{*}$ is, as expected, one of the strong singletons of $u$. The remarkable fact is that we cannot prove that $u^{*} \equiv v^{*}$ entails $u^{*} =_{E} v^{*}$. This is due to the fact that nothing in the theory can assure that {\it that\/} $m$-atom that belongs to $u^{*}$ {\it is the same\/} $m$-atom that belongs to $v^{*}$, since neither the expression $u = v$ nor $u =_{E} v$ are well formed formulas. Furthermore, it is interesting to recall that the usual Extensionality Axiom, which could be used for expressing this fact, is not an axiom of our theory but, instead, we have the ``weak" axiom {\bf Q26} which talks about indistinguishability only, but not about identity. The impossibility of proving the mentioned result should be not regarded as a deficiency of the theory, but rather as expressing that it is closer to what happens in quantum physics than usual set theories. We shall be back to this point below. 

\subsection{Permutations are not observable}

	The next theorem states in the theory ${\cal Q}$ the intuitive idea mentioned above that {\em permutations are not observable\/}. To understand the meaning of this, recall that in standard set theories if $z\in x$, then $(x-\{ z\})\cup \{ w\} = x$ iff $z=w$. So, let us prove the next theorem:

\begin{teorema}
Let $x$ be a qset such that $x \not=_{E} [z]$ and $z$ an $m$-atom such that $z \in x$. If $w \equiv z$ and $w \notin x$, then there exists $w'$ such that $$(x - z') \cup w' \equiv x$$
\end{teorema}

\begin{description}
\item [Proof:] Case 1: $t \in z'$ does not belong to $x$. In this case, $x - z' =_{E} x$ and so we may admit the existence of $w'$ such that its unique element $s$ does belong to $x$ (for instance, $s$ may be $z$ itself); then $(x - z') \cup w' =_{E} x$. 

\noindent
Case 2: $t \in z'$ does belong to $x$. Then $qc(x - z') =_{E} qc(x) - 1$ by the above Theorem. Then we take $w'$ such that its  element is $w$ itself, and so it results that $(x - z') \cap w' =_{E} \emptyset$. Hence, by {\bf Q25},  $qc((x - z') \cup w') =_{E} qc(x)$. This intuitively says that both $(x - z') \cup w'$ and $x$ have the same quantity of indistinguishable elements So, by applying {\bf Q27} (see above), we obtain the theorem. $\Box$
\end{description}

	When $w \notin x$, we have the desired case according to which the theorem is intuitively saying that we have `exchanged' an element of $x$ by an indistinguishable one, and that the resulting fact is that `nothing has occurred at all'. In other words, the resulting qset is indistinguishable from the original one.  The above theorem is the quasi-set theoretical version of the quantum mechanical fact which expresses that permutations of indistinguishable particles are not regarded as observable, as expressed by the so called Indistinguishability Postulate. The relations between quasi-sets and quantum objects are discussed from different points of view in \cite{Krause-95,French-99,Krause-99}. 

\subsection{The Axiom of Choice}

	Finally, the theory ${\cal Q}$ has a version of the axiom of choice.\\

\noindent
{\bf (Q28)} The Axiom of Choice

$\forall_{Q} x (E(x) \wedge \forall y \forall z
(y \in x \wedge z \in x \Rightarrow y \cap z =_{E}  \emptyset  \wedge y 
\neq_{E} \emptyset) \Rightarrow$

$\exists_{Q} u \forall y \forall v (y \in x \wedge v \in y \Rightarrow
\exists_{Q} w (w \subseteq [v] \wedge qc(w) =_{E} 1 \wedge w \cap y 
\equiv w \cap u)))$\\

	Of course this axiom is formulated only to keep ${\cal Q}$ strong enough to be compared with standard ZF, as we did with the Replacement Axioms. As we see, in the axiom, the ``choice qset'' is formed by taking one indistinguishable from each member of the qset $x$. Since we can obtain qsets with quasi-cardinal $2$ whose elements are indistinguishable $m$-atoms, we may reason as if these qsets act as Fraenkel's ``cells'' \cite{Fraenkel-22} in order to obtain, as he did, a proof of the independence of the axiom of choice from the remaining axioms of ${\cal Q}$. As it is well known, the {\em Urelemente\/} of ZFU set theory are indistinguishable in a sense (they are invariant under automorphisms), but even so they do obey the classical theory of identity (they are individuals in a sense). By the contrary, $m$-atoms act as ``legitimate'' indistinguishable entities, so being closer to intuition.

\section{Some Issues of Identity}

	In this section we make a critical analysis on the foundations of equality in first order theories from the point of view of quasi-sets. This is done by means of the study of alternative formulations of quasi-set theory. According to axiom {\bf Q4}:

$$\forall x \forall y (x =_{E} y \Rightarrow (A(x,x) \Rightarrow
A(x,y))).$$

	The question we want raise is: what would happen if we rephrase this sentence in a somehow stronger way? One possibility is:

$$\forall x \forall y (x\equiv y \Rightarrow (A(x,x) \Rightarrow
A(x,y))), \mbox{ with the usual restrictions}.$$

	But in this case, we can easily see that indistinguishability collapses to identity. So, is there any other possible alternative version for {\bf Q4}, which is stronger than {\bf Q4}, but such that it does not collapse to identity. We do not think so. We will illustrate our ideas by means of an alternative version for {\bf Q4}, which we call {\bf Q4$\#$}:

\begin{description}
\item [Q4$\#$ - ] $\forall x \forall y (\neg m(x)\wedge \neg m(y) \wedge x \equiv y \Rightarrow (A(x,x) \Rightarrow
A(x,y))).$
\end{description}

	We can see that {\bf Q4}$\#$ is stronger than {\bf Q4} since {\bf Q4}$\#$ allows substitutivity for indistinguishable qsets which are not extensionaly identical. But the question is: is {\bf Q4}$\#$ weaker than identity?

	In this section we will prove some lemmas and theorems within the scope of a quasi set theory which replaces {\bf Q4} in ${\cal Q}$ by {\bf Q4}$\#$. We call this theory ${\cal Q}\#$. Our main goal is to allow a better understanding of identity in first order theories. We will prove here that ${\cal Q}\#$ is equivalent to ZFU set theory, since indistinguishability, in this case, collapses to identity.

\begin{lema}\label{cachaca}
For all qsets $x$ and $y$ we have:

\begin{enumerate}

\item If $t\in x$ and $x\equiv y$ then $t\in y$;

\item If $x\equiv y$ then $x=_E y.$

\end{enumerate}
\end{lema}

\begin{description}
\item [Proof:] (1) Since $x$ is a qset, then $\neg m(x)$, according to definition (\ref{defqset}). Since $x\equiv y$ then $\neg m(y)$ ({\bf Q9}). Now, for all qsets $u$ and $v$ let $A_t(u,v) := (t\in u\wedge t\in v)$. By hypothesis, $t\in x \wedge t\in x$, i.e., $A_t(x,x)$. Since $x\equiv y$ then, according to {\bf Q4$\#$}, $A_t(x,y)$. Hence $t\in y$. (2) If $x\equiv y$ then, according to item 1 of this proof, $t\in x$ iff $t\in y$, which means that $x=_Ey$.$\Box$
\end{description}

\begin{lema}\label{aguardente}
For all $m$-atom $x$, $[x]=_Ex^*$, where $x^*$ is given by definition (\ref{star}).
\end{lema}

\begin{description}
\item [Proof:] Let $t\in [x]$. According to the definition of weak singleton, $t\equiv x$. According to {\bf Q9}, $t$ is an $m$-atom. From lemma (\ref{pinga}), we have $t^*\equiv x^*$. Note that the proof of lemma (\ref{pinga}) does not make any reference to {\bf Q4}, so we can use it here, although we are working in ${\cal Q}\#$. From definition (\ref{star}) we know that $t^*$ and $x^*$ are qsets. So, $t^*=_Ex^*$, according to lemma (\ref{cachaca}). So, since $t\in t^*$, then $t\in x^*$. Therefore $[x]\subseteq x^*$. From lemma (\ref{Lemma}) (which also does not make any reference to {\bf Q4} in its proof) $x^*\subseteq [x]$. So, $[x]=_E x^*$.$\Box$.
\end{description}

\begin{lema}\label{aguardente2}
For all $m$-atoms $x$ and $y$, the following conditions are equivalent:

\begin{enumerate}

\item $x\equiv y$;

\item $x^*\equiv y^*$;

\item $x^* =_E y^*$;

\item $[x] =_E [y]$;

\item $[x] \equiv [y]$.

\end{enumerate}
\end{lema}

\begin{description}
\item [Proof:] By lemma (\ref{pinga}) (1)$\Rightarrow$ (2); By {\bf Q4}$\#$ (2)$\Rightarrow$(3); By lema (\ref{aguardente}) (3)$\Rightarrow$ (4); by theorem (6), item (3), (4)$\Rightarrow$ (5); By the same theorem, item (4), (5)$\Rightarrow$ (1).$\Box$
\end{description}

\begin{teorema}\label{churrasco}
For all $m$-atom $x$, $qc([x]) =_E 1$.
\end{teorema}

\begin{description}
\item [Proof:] Straightforward from lemma (\ref{aguardente}) and theorem (\ref{theorem}) (which makes no reference to {\bf Q4} in its proof).$\Box$
\end{description}

\begin{lema}\label{lema1}
For all $m$ atoms $x$ and $y$ and for all qset $w$, if $x\equiv y$ and $x\in w$, then $y\in w$.
\end{lema}

\begin{description}
\item [Proof:] From lemma (\ref{pinga}) $x^*\subseteq w$. From lemma (\ref{aguardente2}) $y^*\subseteq w$. From lemma (\ref{Lemma}), item (4), $y\in w$.$\Box$
\end{description}

\begin{lema}\label{lema2}
For all $m$-atoms $x$ and $y$ and for all qset $w$:

\begin{enumerate}

\item If $x\equiv y$ and $[x]\in w$, then $[y]\in w$;

\item If $x\equiv y$ and $[x]\subseteq w$ then $[y]\subseteq w$.

\end{enumerate}
\end{lema}

\begin{description}
\item [Proof:] If $x\equiv y$ then, from lemma (\ref{aguardente2}), $[x]\equiv [y]$. But $[x]$ and $[y]$ are qsets. Then, from {\bf Q4}$\#$ we have (1) and (2).$\Box$
\end{description}

\begin{lema}\label{lema3}
For all $m$-atoms $x$ and $y$ and for all $\lambda$ and $z$:

\begin{enumerate}

\item If $x\equiv y$ and $qc([x])=_E\lambda$ then $qc([y])=_E\lambda$;

\item If $x\equiv y$ and $z\in [x]$ then $z\in [y]$.

\end{enumerate}
\end{lema}

\begin{description}
\item [Proof:] (1) Follows  from theorem (\ref{churrasco}). (2) Follows from {\bf Q12} and {\bf Q2}.$\Box$
\end{description}

	Next follows the main result of this section:

\begin{teorema}
With the usual restrictions, $\forall x\forall y(x\equiv y \Rightarrow (A(x,x)\Rightarrow A(x,y)))$.
\end{teorema}

\begin{description}
\item [Proof:] Suppose $m(x)$. From {\bf Q9}, $m(y)$. In this case, $A(x,x)$ is only built from the following types of atomic formulas for some qset $w$: (1) $x\in w$; (2) $[x]\in w$; (3) $[x]\subseteq w$; (4) $qc([x])=_E\lambda$; (5) $z\in [x]$. From lemmas (\ref{lema1}), (\ref{lema2}), and (\ref{lema3}), we have: (1) $y\in w$; (2) $[y]\in w$; (3) $[y]\subseteq w$; (4) $qc([y])=_E\lambda$; (5) $z\in [y]$, i.e., if $A(x,x)$ then $A(x,y)$. Suppose now that $\neg m(x)$. According to {\bf Q9} $\neg m(y)$. By means of {\bf Q4}$\#$, if $A(x,x)$ then $A(x,y)$.$\Box$
\end{description}

	This last theorem says that indistinguishability $\equiv$ collapses into identity in ${\cal Q}\#$. So, ${\cal Q}\#$ is equivalent to standard ZFU.

	This last theorem depends essentially on the Weak Axiom of Extensionality {\bf Q26}. It is worth to remark that axiom {\bf Q4}$\#$ was used by one of us in \cite{Krause-92}. But in that paper the Axiom of Extensionality was different also. So, some results presented here are not valid in the quasi-set theory introduced in \cite{Krause-92}.

\section{Physics: The Maxwell-Boltzmann Statistics}

	The contents of this section are also discussed in \cite{Sant'Anna-00b}.

	According to usual textbooks on statistical mechanics, Maxwell-Boltzmann (MB) statistics gives us the most probable distribution of $N$ {\em distinguishable\/} objects into, say, boxes with a specified number of objects in each box. In this section we show that the hypothesis concerning the objects being distinguishable is unnecessary.

\subsection{Some Standard Results in ZF}

	It is a well known theorem in Zermelo-Fraenkel set theory the following:

\begin{lema}
If $x$ is a finite ZF-set, then $$card({\cal P}(x)) = 2^{card(x)}.$$\label{card2}
\end{lema}

\begin{teorema}
Let $x$ be a non-empty and finite ZF-set. If we define $x_2$ as a set of ordered pairs $\langle y_1,y_2\rangle$ such that $y_1,y_2\in {\cal P}(x)$, $y_1\cup y_2 = x$, and $y_1\cap y_2 = \emptyset$ then $card(x_2) = 2^{card(x)}$.\label{yyy}
\end{teorema}

	This theorem corresponds to say that the number of ways we can distribute $N$ distinguishable particles ($N = card(x)$) between {\em two\/} boxes (represented by the ordered {\em pair\/} $\langle y_1,y_2\rangle$) is $2^N$.

\begin{teorema}
Let $x$ be a finite ZF-set such that $card(x) = N$. If we define $x_n$ as a set of ordered $n$-tuples $\langle y_1,\cdots, y_n\rangle$ such that for all $i = 1,\cdots, n$ we have $y_i\in{\cal P}(x)$, $\bigcup_i y_i = x$, and $i\neq j\Rightarrow y_i\cap y_j = \emptyset$, then $card(x_n) = n^N$.\label{nN}
\end{teorema}

	We could rewrite theorem (\ref{nN}) as:

\begin{teorema}
Let $x$ be a finite ZF-set such that $card(x) = N$. If we define $x_n$ as a set of ordered $n$-tuples $\langle y_1,\cdots, y_n\rangle$ such that for all $i = 1,\cdots, n$ we have $y_i\in{\cal P}(x)$, $\bigcup_i y_i = x$, and $\sum_i card(y_i) = card(x)$, then $card(x_n) = n^N$.\label{nNp}
\end{teorema}

\begin{description}
\item [Proof:] Analogous to the proof of theorem (\ref{nN}), since $\bigcup_i y_i = x$, and $i\neq j\Rightarrow y_i\cap y_j = \emptyset$ iff $\bigcup_i y_i = x$, and $\sum_i card(y_i) = card(x)$.$\Box$
\end{description}

	This theorem corresponds to say that the number of ways that we can distribute $N$ distinguishable particles ($N = card(x)$) among $n$ boxes (represented by the ordered $n$-tuple $\langle y_1,\cdots , y_n\rangle$) is $n^N$.

\subsection{Quasi-Set-Theoretical Combinatorics}

	We can obtain a perhaps more fruitful theory, which allows us a quasi-set theoretical combinatorics, if we exchange the axiom {\bf Q25} by the following postulate, which is a generalization of {\bf Q25}), as well as a quasi-set theoretical version of theorem (\ref{nN}):

\begin{description}

\item[Q25'] Let $x$ be a finite quasi-set such that $qc(x) =_E N$. If we define $z_n$ as the quasi-set whose elements are ordered $n$-tuples $\langle y_1,\cdots , y_n\rangle$, where, for all $i =_E 1,\cdots, n$, we have $y_i\in {\cal P}(x)$, $\bigcup_{i }y_i = x$, and $\sum_i qc(y_i) =_E qc(x)$, then we have the following:

\begin{equation}
qc(z_n) =_E n^N.\label{MBgeneral}
\end{equation}

\end{description}

	In the case where $n =_E 2$, we have a sentence which is equivalent to axiom {\bf Q25}.

	The main role of axiom {\bf Q25'} is to allow us a quasi-set theoretical combinatorics which can be useful to cope with distribution functions. From the mathematical point of view, it is important to show that the substitution of axiom {\bf Q25} by axiom {\bf Q25'} does not entail any inconsistency in quasi-set theory. This shall be proved in the Section 5. The point, at this moment, is that {\bf Q25} is very `poor' if we are interested on a quasi-set-theoretical combinatorics with more than two physical states or `boxes', as exemplified in the Introduction. Besides, axiom {\bf Q25'} is our quasi-set theoretical version of theorem (\ref{nNp}).

	If we recall the polynomial of Leibniz, we can rewrite equation (\ref{MBgeneral}) as:

\begin{equation}
qc(z_n) =_E n^N =_E \sum\frac{N!}{\Pi_{i = 1,\cdots n}n_i!},\label{MBgeneral2}
\end{equation}

\noindent
where the sum is over all possible combinations of nonnegative integers $n_i$ such that $\sum_{i = 1,\cdots, n}n_i =_E N$.

	If we interpret $n$ as the number of physical states, $N$ as the total number of particles and $n_i$ as the number of particles associated to each physical state $i$, then it is easy to see that each parcel of the summation in equation (\ref{MBgeneral2}) is a possible MB distribution of $N$ particles among $n$ states. The most probable among all these parcels is the MB distribution. So, we can add equation (\ref{MBgeneral2}), with its respective interpretation, as another extra-assumption in quasi-set theory. In other words, we are generalizing theory ${\cal Q}$, by replacing axiom {\bf Q25} by axiom {\bf Q25'}. We refer to this generalized quasi-set theory as ${\cal Q'}$. If we do not replace axiom {\bf Q25} by {\bf Q25'}, there is no manner of saying anything about a distribution of $N$ particles among an arbitrary number $n$ of states or boxes. In this case, we would be confined to the very particular case of 2 states.

	It is easy to see that, for all $i$ we have $n_i =_E qc(y_i)$. Axiom {\bf Q25'} is just another way of saying that the number of ways we can distribute $N$ objects (either distinguishable or not) among $n$ boxes is $n^N$. The condition that $\bigcup_{i }y_i =_E x$, and $\sum_i qc(y_i) =_E qc(x)$ is simply a way to guarantee that there will be no `repeated occurence' of the same object in two boxes. Nevertheless, it is obvious that the expression `repeated occurence', in this quasi-set-theoretical context, is just an intuitive approach for didactical purposes, since there is no sense in saying that the `same' object cannot occupy two boxes.

	The reader could ask: what are the so-called ``boxes''? Each $y_i$ corresponds to a given box or physical state. There can be, of course, two indistinguishable boxes $y_i$ and $y_j$. In this case, the labels $i$ and $j$ cannot individualize each box. They are just different names, or labels, attributed to two indistinguishable objects (qsets, in this case).

\subsection{One Simple Example}

	Now, let us exhibit an example in order to illustrate our ideas. Consider a collection of three indistinguishable particles to be distributed between two possible states or `boxes'. According to standard textbooks on statistical mechanics, there are only four possibilities of distribution. On the other hand, according to our axiomatic framework -- axiom {\bf Q25'} -- there are eight possibilities. If we impose that the occupation number of each box is constant, the number of possibilities corresponds to one parcel of the sum in equation (\ref{MBgeneral2}).

	The question now is: what about the extra four possibilities predicted by axiom {\bf Q25'}? The eight possibilities predicted by {\bf Q25'} and equation (\ref{MBgeneral2}) come from 

$$2^3 = \frac{3!}{3!0!} + \frac{3!}{2!1!} + \frac{3!}{1!2!} + \frac{3!}{0!3!}.$$

	So, we have one possibility with 3 particles in the first state and no particle in the second state, plus three {\em indistinguishable\/}  possibilities with 2 particles in the first state and 1 particle in the second state, plus three {\em indistinguishable\/} possibilities with 1 particle in the first state and 2 particles in the second state, plus one single possibility with no particle in the first state and 3 particles in the remaining one. The calculation of the most probable case is made for a large number of particles, following the standard calculations of statistical mechanics.

	Following our example, axiom {\bf Q25'} says that we can distribute $3$ objects (either indistinguishable or not) among $2$ boxes in $2^3$ manners (either indistinguishable or not). But this axiom does not say {\em how\/} can we make this distribution. If we do not appeal to equation (\ref{MBgeneral2}), we have the following: according to Fig. 1, there are, at least, by means of axiom {\bf Q16}, {\em four\/} possible distributions. But axiom {\bf Q25'} says that there are eight possible distributions. One possibility is something like Fig. 2, that is, the four distributions in Fig. 1 {\em plus\/} four distributions which are indistinguishable from the third distribution of Fig. 1. The reader can easily imagine other possibilities. So, axiom {\bf Q25'} by itself does not allow us to derive MB statistics. It simply says that MB statistics is a possibility even in a collection of indiscernibles. Axiom {\bf Q25'} {\em and\/} equation (\ref{MBgeneral2}), with its respective interpretation in the context of {\bf Q25'}, is a way to say that the only possibility is that one illustrated at the Fig. 3.

\section{Quantum Statistics}

	Since we may have MB distribution among non-individuals, what is the difference between quantum statistics and MB, after all? In Bose-Einstein we take into account {\em only distinguishable possibilities\/}, among all possibilities predicted by axiom {\bf Q25'}. And Fermi-Dirac is derived in the same vein, but with the additional assumption of the quasi-set theoretical version of Pauli's Exclusion Principle: $qc(y_i)\leq 1$ for each $i$ in {\bf Q25'}. Put it in another way, quantum statistics may be seen as special cases of MB statistics in a collection of indistinguishable particles. Another way, more complicated, to get quantum statistics in quasi-set theory is introduced in \cite{Krause-99}.

\section{Acknowledgments}

	We acknowledge with thanks the insightful discussions that we had with Jos\'e Renato Ramos Barbosa within the context of the Analice Gebauer Volkov Seminars at Federal University of Paran\'a.

\newpage

\begin{figure}
\renewcommand{\arraystretch}{0.7}
\[
\begin{array}{|c|c|}\hline
\bullet\bullet\bullet & \;\; \\ \hline
\bullet\bullet & \bullet \\ \hline
\bullet & \bullet\bullet \\ \hline
\;\; & \bullet\bullet\bullet \\ \hline
\end{array}
\]
\caption{The `first' four possible distributions of 3 objects (indistinguishable or not) among 2 boxes. Each line represents one possible distribution and each bullet represents an object.}
\end{figure}
\begin{figure}
\renewcommand{\arraystretch}{0.7}
\[
\begin{array}{|c|c|}\hline
\bullet\bullet\bullet & \;\; \\ \hline
\bullet\bullet & \bullet \\ \hline
\bullet & \bullet\bullet \\ \hline
\;\; & \bullet\bullet\bullet \\ \hline
\bullet & \bullet\bullet \\ \hline
\bullet & \bullet\bullet \\ \hline
\bullet & \bullet\bullet \\ \hline
\bullet & \bullet\bullet \\ \hline
\end{array}
\]
\caption{One possible sequence of the eight possible distributions of 3 objects among 2 boxes according to axiom {\bf Q25'}.}
\end{figure}
\begin{figure}
\renewcommand{\arraystretch}{0.7}
\[
\begin{array}{|c|c|}\hline
\bullet\bullet\bullet & \;\; \\ \hline
\bullet\bullet & \bullet \\ \hline
\bullet\bullet & \bullet \\ \hline
\bullet\bullet & \bullet \\ \hline
\bullet & \bullet\bullet \\ \hline
\bullet & \bullet\bullet \\ \hline
\bullet & \bullet\bullet \\ \hline
\;\; & \bullet\bullet\bullet \\ \hline
\end{array}
\]
\caption{The only possible distribution of 3 objects among 2 boxes, if we conjugate axiom {\bf Q25'} and equation (2).}
\end{figure}


\begin{thebibliography}{9}

\bibitem{Browder-76} Browder, F. E. (ed.), {\em Proceedings of the Symposium on Pure Mathematics of the American Mathematical Society - Mathematical Developments Arising from Hilbert Problems)\/} {\bf 28} (AMS, Providence, 1976).

\bibitem{Cantor-55} Cantor, G., {\em Contributions to the Founding of Transfinite Numbers\/} (Dover, New York, 1955).

\bibitem{daCosta-99} da Costa, N. C. A. and Krause, D., `Set theoretical models for quantum systems', in Dalla Chiara, M. L. {\t et al.\/} (eds.), {\it Language, quantum, music\/}, Kluwer Ac. Press, 1999, 171--181. 

\bibitem{Enderton-77} Enderton, H. B., {\it Elements of set theory\/}, Academic Pres, 1977.

\bibitem{Fraenkel-22} Fraenkel, A. A., `The notion of `definite' and the independence of the axiom of choice' (1922), in J. van Heijennort (ed.) {\em From Frege to G\"odel: a source book in mathematical logic, 1879--1931\/} (Harvard Un. Press, 1967) 284--289.

\bibitem{French-95}  French, S. and Krause, D., `Vague identity and quan\-tum  non-indi\-vidual\-ity', {\em Analysis\/} {\bf 55} (1), 1995, 20--26.

\bibitem{French-99}  French, S. and Krause, D., `The logic of quanta', in Cao, T. Y. (ed.), {\it Conceptual foundations of quantum field theory\/}, Cambridge Un. Press, 1999, 324--242.

\bibitem{French-9*} French, S., Krause, D. and Maidens, A., `Quantum vagueness', preprint, University of Leeds.

\bibitem{Heisenberg-89} Heisenberg, W., `What is an elementary particle?', in Heisenberg, W., {\it Encounters with Einstein and other essays on people, places, and particles\/}, Princeton Un. Press, 1989.

\bibitem{Huggett-99} Huggett, N., `Atomic metaphysics', {\em The Journal of Philosophy\/} {\bf 96} 5-24 (1999).

\bibitem{Krause-90} Krause, D., {\it N\~ao-reflexividade, indisting\"uibilidade e agregados de Weyl\/}, Thesis, FFLCH-USP, 1990.

\bibitem{Krause-92} Krause, D.,`On a  quasi-set theory', {\em Notre Dame Journal of Formal Logic\/} {\bf 33} (3), 1992, 402-411.

\bibitem{Krause-96a} Krause, D., `Axioms for collections of indistinguishable objects',  {\it Logique et Analyse\/} {\bf 153--154}, 1996, 69--93.

\bibitem{Krause-95} Krause, D.  and French, S., `A formal fra\-me\-work
\- for \- quan\-tum non--\-indi\-vidual\-ity', {\em Synthese\/} {\bf 102},
1995, 195--214. 

\bibitem{Krause-96b} Krause, D. and French, S., `Quantum objects are vague objects', {\it Sorites\/} {\bf 6}, 1996, 21--33. 

\bibitem{Krause-99} Krause, D., Sant'Anna, A. S. and Volkov, A. G., `Quasi-set theory for bosons and fermions', {\it Found. Phys. Lett.\/}, {\bf 12} (1), 1999, 51--66. 

\bibitem{Manin-76} Manin, Yu. I., `Problems of Present Day Mathematics I: Foundations', in Browder, F. E. (ed.), {\em Mathematical problems arising from Hilbert problems\/}, Proceedings of Symposia in Pure Mathematics XXVIII, Providence, AMS, 1976, p. 36.

\bibitem{Manin-77} Manin, Yu. I., {\it A course in mathematical logic\/}, New York, Springer-Verlag, 1977. 

\bibitem{Mendelson-97} Mendelson, E., {\it Introduction to mathematical logic\/}, London, Chapman \& Hall, 4th. ed., 1997.

\bibitem{Post-63} Post, H., `Individuality and physics', {\em The Listener\/} {\bf 70} 534-537 (1963).

\bibitem{Sant'Anna-00a} Sant'Anna, A. S., `Elementary particles, hidden variables, and hidden predicates', {\em Synthese\/} {\bf 125} 233-245 (2000).

\bibitem{Sant'Anna-00b} Sant'Anna, A. S. and A. M. S. Santos, `Quasi-set-theoretical foundations of statistical mechanics: a research program', {\em Found. Phys.\/}, {\bf 30} 101-120 (2000).

\bibitem{Schrodinger-52} Schr\"odinger, E., {\em Science and humanism\/}, Cambridge Un. Press, Cambridge, 1952. 

\bibitem{Schrodinger-98} Schr\"odinger, E., `What is an elementary particle?', (reprinted in) Castellani, E. (ed.), {\it Interpreting bodies: classical and quantum objects in modern physics\/}, Princeton, Princeton Un. Press, 1998, 197--210.

\bibitem{vanFraassen-91} van Fraassen, B. C., {\it Quantum mechanics: An empiricist view\/},  Oxford, Clarendon Press, 1991. 

\bibitem{Weingartner-96} Weingartner, P., `Under what transformations are laws invariants?', in P. Weingartner and G. Schurz (eds.) {\em Law and Prediction in the Light of Chaos Research\/} Lecture Notes in Physics 473 (Springer, New York, 1996) 47-88.

\bibitem{Weyl-49} Weyl, H., {\it Philosophy of mathematics and natural science\/}, Princeton Un. Press, 1949. 

\end{thebibliography}
\end{document}